\documentclass{maxart-noswap}
\usepackage{amssymb,amsfonts,amsxtra}

\theoremstyle{plain}
\newtheorem{proposition}{Proposition}
\newtheorem{theorem}[proposition]{Theorem}

\newtheorem{lemma}[proposition]{Lemma}

\theoremstyle{definition}
\newtheorem{definition}[proposition]{Definition}
\newtheorem{remark}[proposition]{Remark}
\newtheorem*{ack}{Acknowledgements}
\newtheorem*{example}{Main Example}

\numberwithin{proposition}{section}
\renewcommand{\theenumi}{\alph{enumi}}

\newcommand{\zindex}[3]{\put(#1,#2){\makebox(0,0){${#3}$}}}
\newenvironment{pict}[2]%
	{\setlength{\unitlength}{1mm}
	\begin{center}
	\begin{picture}(#1,#2)
	\scriptsize}{\end{picture}
	\end{center}

	\noindent}

\DeclareFontFamily{OML}{rsfs}{\skewchar\font'177}
\DeclareFontShape{OML}{rsfs}{m}{n}{ <5> <6> rsfs5 <7> <8> <9> rsfs7
  <10> <10.95> <12> <14.4> <17.28> <20.74> <24.88> rsfs10 }{}
\DeclareMathAlphabet{\mathfs}{OML}{rsfs}{m}{n}

\newcommand{\As}{{\mathfs A}}
\newcommand{\Cs}{{\mathfs C}}
\newcommand{\Gs}{{\mathfs G}}
\newcommand{\R}{{\mathbb R}}
\newcommand{\Z}{{\mathbb Z}}
\renewcommand{\geq}{\geqslant}

\renewcommand{\today}{April 10, 2003}

\begin{document}

\renewcommand{\thefootnote}{\fnsymbol{footnote}}
\title[On uniqueness of JSJ decompositions]{On uniqueness of JSJ
decompositions \\ of finitely generated groups} 
\author{Max Forester\footnote{Research supported by EPSRC grant
GR/N20867}} 
\address{Mathematics Institute, University of Warwick, Coventry, CV4 7AL,
UK} 
\email{forester@maths.warwick.ac.uk}

\begin{abstract}
We give an example of two JSJ decompositions of a group that are not
related by conjugation, conjugation of edge--inclusions, and slide
moves. This answers the question of Rips and Sela stated in
\cite{ripssela}. 

On the other hand we observe that any two JSJ decompositions of a group
are related by an elementary deformation, and that strongly
slide--free JSJ decompositions are genuinely unique. These results hold
for the decompositions of Rips and Sela, Dunwoody and Sageev, and
Fujiwara and Papasoglu, and also for accessible decompositions. 
\end{abstract}

\maketitle 
\renewcommand{\theauthors}{Max Forester}

\section*{Introduction}

In this paper we discuss the extent to which JSJ decompositions of
groups are unique. In \cite{sela:hypjsj2} and \cite{ripssela} it was
shown that if $G$ is a (Gromov) hyperbolic group then any two JSJ
decompositions of $G$ must be related by conjugation, conjugation of
edge--inclusions, and slide moves. Rips and Sela also noted that the same
uniqueness statement holds in many other cases. However, the general case
was left open as a question. Here we give examples of JSJ decompositions
of a finitely presented group that are not related by such moves,
answering their question in the negative. 

In light of the examples it is natural to ask what form of uniqueness
does hold for finitely presented groups. It turns out that JSJ
decompositions are unique up to \emph{elementary deformation}, a notion
that is studied extensively in \cite{forester:trees}. Furthermore, if a
decomposition is \emph{strongly slide--free} then it is genuinely unique. As
we will see in Section \ref{uniquesec}, these results follow directly
from the more general results in \cite{forester:trees}. The same
uniqueness results also hold for accessible (or one-ended)
decompositions. 

In this paper we focus mainly on the JSJ decomposition of Rips and
Sela, though our results apply equally well to the JSJ decompositions of
Dunwoody and Sageev \cite{dunwoodysageev} and Fujiwara and Papasoglu
\cite{fujiwara:jsj}. It seems worthwhile to mention the uniqueness
properties of other related decompositions. 

There have been several constructions, originally inspired by the
canonical decompositions of $3$--manifolds due to Jaco and Shalen, and
Johannson \cite{jacoshalen,johannson}. The first result of this kind was
Kropholler's decomposition for Poincar\'e duality groups in
\cite{krop:decomposition}. This decomposition is unique. Sela defined a
JSJ decomposition for torsion free hyperbolic groups in
\cite{sela:hypjsj2}, and Bowditch found an equivalent topological
construction for one-ended hyperbolic groups; this decomposition is
essentially unique \cite{bowditch:cutpoints,srilatha}. Next came the
three decompositions already mentioned
\cite{ripssela,dunwoodysageev,fujiwara:jsj}, whose uniqueness properties
are discussed in this paper. Finally there is a recent construction due to
Scott and Swarup, described in the $3$--manifold case in
\cite{scottswarup}, and in the general case in
\cite{scottswarup:decomp}. This decomposition is unique and it agrees
with the topological JSJ decomposition in the case of a $3$--manifold
group. On the other hand, the JSJ decompositions considered here are
sometimes finer, and can reveal more information about the group. 

The paper is organized as follows. In Section \ref{prelimsec} we discuss
moves between decompositions. In Section \ref{examplesec} we present the
main examples, which are generalized Baumslag--Solitar trees. Most of the
section is devoted to showing that such trees qualify as JSJ
decompositions, under mild assumptions (Theorem \ref{gbsjsj} and
Proposition \ref{gbsunfolded}). This result is interesting in its own
right. In Section \ref{uniquesec} we discuss our uniqueness result,
Theorem \ref{uniqueness}. 

\begin{ack}
I would like to thank Gilbert Levitt for pointing out that the results of
\cite{forester:trees} could be applied to JSJ decompositions. Similar
suggestions were made by M. Sageev, Z. Sela, and G. A. Swarup. During
this work I benefitted from discussions and correspondence with Peter
Scott and G. A. Swarup.  I also thank David Epstein for his
encouragement, and the referee for suggesting improvements to the
exposition. 
\end{ack}

\section{Preliminaries} \label{prelimsec} 

We will use Serre's notation for graphs and trees. Thus a graph $A$ is a
pair of sets $(V(A)$, $E(A))$ with maps $\partial_0, \partial_1 \co E(A)
\to V(A)$ and $e \mapsto \overline{e}$ (for $e \in E(A)$), such that
$\partial_i\overline{e} = \partial_{1-i} e$ and $e \not= \overline{e}$
for all $e$. An element $e \in E(A)$ is to be thought of as an oriented
edge with initial vertex $\partial_0 e$ and terminal vertex $\partial_1
e$. 

Let $G$ be a group. A \emph{$G$--tree} is a tree with a
$G$--action by automorphisms, without inversions. There is a
correspondence between $G$--trees and graphs of groups having
fundamental group $G$, as explained in \cite{serre:trees}. We will
consider certain moves between 
graphs of groups that do not change the fundamental group. Equivalently,
these are moves between $G$--trees. A more
complete discussion of these moves is given in \cite[\S 3]{forester:trees}. 

\begin{definition}
In a \emph{collapse move}, an edge in a graph of
groups carrying an amalgamation of the form $A \ast_C C$ is collapsed to
a vertex with group $A$. Every inclusion map having target group $C$ is
reinterpreted as a map into $A$, via the injective map of vertex groups 
$C \hookrightarrow A$. 
This move simplifies the underlying graph
without enlarging any vertex or edge groups. 

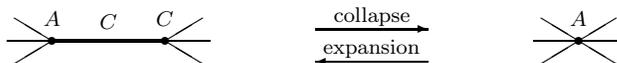
\begin{pict}{90}{10}
\thicklines
\put(10,5){\circle*{1}}
\put(25,5){\circle*{1}}

\put(10,5){\line(1,0){15}}

\thinlines
\put(45,6.5){\vector(1,0){15}}
\put(60,2.2){\vector(-1,0){15}}
\zindex{52.5}{8.2}{\mbox{collapse}}
\zindex{52.5}{4}{\mbox{expansion}}

% \put(25,5){\line(3,5){3}}
\put(25,5){\line(5,-3){5}}
\put(25,5){\line(5,3){5}}

\put(10,5){\line(-5,3){5}}
\put(10,5){\line(-5,-3){5}}

\put(80,5){\circle*{1}}
\put(80,5){\line(-5,3){5}}
\put(80,5){\line(-5,-3){5}}

% \put(80,5){\line(3,5){3}}
\put(80,5){\line(5,-3){5}}
\put(80,5){\line(5,3){5}}

\put(4,5){\line(1,0){27}}
\put(74,5){\line(1,0){12}}

\zindex{10}{8}{A}
\zindex{17.5}{7.5}{C}
\zindex{25}{8}{C}

\zindex{80}{8}{A}
\end{pict} 
An \emph{expansion move} is the reverse of a collapse move. Both of these
moves are called \emph{elementary moves}. An \emph{elementary
deformation} is a finite sequence of such moves. A graph of groups is
\emph{reduced} if it admits no collapse moves. This means that if an
inclusion map from an edge group to a vertex group is an isomorphism,
then the edge is a loop. 
\end{definition}

\begin{definition}
The elementary deformation  shown below, consisting of an expansion move
followed by a collapse, is called a \emph{slide move}. In order to
perform the expansion it is required that $D \subseteq C$ (regarded as
subgroups of $A$). 

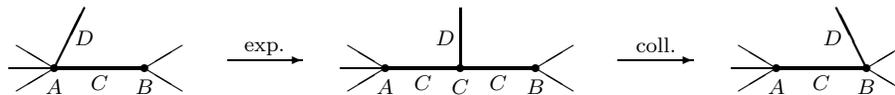
\begin{pict}{120}{12}
\thicklines
\put(102,4){\circle*{1}}
\put(114,4){\circle*{1}}
\put(102,4){\line(1,0){12}}
\put(114,4){\line(-1,2){4}}

\zindex{102}{1.5}{A}
\zindex{114}{1.5}{B}
\zindex{108}{2}{C}
\zindex{109.5}{8}{D}

\thinlines
\put(114,4){\line(5,3){5}}
\put(114,4){\line(5,-3){5}}
\put(102,4){\line(-5,3){5}}
\put(96,4){\line(1,0){6}}
\put(102,4){\line(-5,-3){5}}

%%%%%%%%

\thicklines
\put(6,4){\circle*{1}}
\put(18,4){\circle*{1}}
\put(6,4){\line(1,0){12}}
\put(6,4){\line(1,2){4}}

\zindex{6}{1.5}{A}
\zindex{18}{1.5}{B}
\zindex{12}{2}{C}
\zindex{10}{8}{D}

\thinlines
\put(18,4){\line(5,3){5}}
\put(18,4){\line(5,-3){5}}
\put(6,4){\line(-5,3){5}}
\put(0,4){\line(1,0){6}}
\put(6,4){\line(-5,-3){5}}

%%%%%%%%

\thicklines
\put(50,4){\circle*{1}}
\put(60,4){\circle*{1}}
\put(70,4){\circle*{1}}
\put(50,4){\line(1,0){20}}
\put(60,4){\line(0,1){8}}

\zindex{50}{1.5}{A}
\zindex{55}{2}{C}
\zindex{60}{1.5}{C}
\zindex{65}{2}{C}
\zindex{70}{1.5}{B}
\zindex{58}{8}{D}

\thinlines
\put(70,4){\line(5,3){5}}
\put(70,4){\line(5,-3){5}}
\put(50,4){\line(-5,3){5}}
\put(44,4){\line(1,0){6}}
\put(50,4){\line(-5,-3){5}}

%%%%%%%%

\zindex{34}{6.5}{\mbox{exp.}}
\zindex{86}{7}{\mbox{coll.}}
\put(29,5){\vector(1,0){10}}
\put(81,5){\vector(1,0){10}}

\end{pict} 
It is permitted for the edge carrying $C$ to be a
loop; in this case the only change to the graph of groups is in the
inclusion map $D \hookrightarrow A$. See Proposition \ref{gbsmoves} for
an example. 
\end{definition}

\begin{definition}
A \emph{fold} is most easily described in terms of $G$--trees. The graph of
groups description involves many different cases which are explained in
\cite{bestvina:accessibility}. To perform a fold
in a $G$--tree one chooses edges $e$ and $f$ with $\partial_0 e =
\partial_0 f$, and identifies $e$ and $f$ to a single edge. One also
identifies $\gamma e$ with $\gamma f$ for every $\gamma \in G$, so that
the resulting quotient graph has a $G$--action. It is not difficult to
show that the new graph is a tree. 
\end{definition}

\begin{definition}
A \emph{generalized Baumslag--Solitar tree} is a $G$--tree whose vertex
and edge stabilizers are all infinite cyclic. The groups $G$ that arise
are called \emph{generalized Baumslag--Solitar groups}. Examples include
the classical Baumslag--Solitar groups and torus knot groups. When
discussing specific examples it is convenient to use edge--indexed
graphs, as seen in the next section. They depict graphs of groups in
which all edge and vertex groups are $\Z$. The indices define the
inclusion maps, which are simply multiplication by various non-zero
integers. 
\end{definition}

\section{Two JSJ decompositions not related by slide moves}
\label{examplesec} 

In \cite[p. 106]{ripssela}, Rips and Sela ask whether any two JSJ
decompositions of a group must be related by conjugation, conjugation of
edge--inclusions, and slide moves. If one regards the JSJ decomposition
as a $G$--tree then the first two modifications have no effect (up to
$G$--isomorphism). Thus, they are asking whether two such $G$--trees are
related by slide moves. We show by example that this need not be the
case. 

The examples are generalized Baumslag--Solitar trees. We will describe
two such trees and show that they are related by an elementary
deformation. This implies that the two groups are the same. Then we
verify that the trees are not slide--equivalent, and that they represent
JSJ decompositions of their common group. 

\begin{proposition} \label{gbsmoves} 
If an elementary move is performed on a generalized
Baumslag--Solitar tree, then the quotient graph of groups changes locally
as follows: 

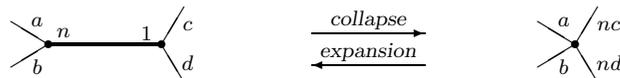
\begin{pict}{100}{10}
\thicklines
\put(10,5){\circle*{1}}
\put(25,5){\circle*{1}}

\put(10,5){\line(1,0){15}}

\thinlines
\put(45,6.5){\vector(1,0){15}}
\put(60,2.2){\vector(-1,0){15}}
\zindex{52.5}{8.2}{\mbox{collapse}}
\zindex{52.5}{4}{\mbox{expansion}}

\put(25,5){\line(3,5){3}}
\put(25,5){\line(3,-5){3}}

\put(10,5){\line(-5,3){5}}
\put(10,5){\line(-5,-3){5}}

\put(80,5){\circle*{1}}
\put(80,5){\line(-5,3){5}}
\put(80,5){\line(-5,-3){5}}

\put(80,5){\line(3,5){3}}
\put(80,5){\line(3,-5){3}}

\scriptsize
\zindex{8.5}{8}{a}
\zindex{8.5}{2}{b}
\zindex{12}{6.5}{n}
\zindex{23}{6.5}{1}
\zindex{28.5}{7.5}{c}
\zindex{28.5}{2.5}{d}

\zindex{78.5}{8}{a}
\zindex{78.5}{2}{b}
\zindex{84.5}{7.5}{nc}
\zindex{84.5}{2.5}{nd}
 
\end{pict} 
A slide move has the following description: 

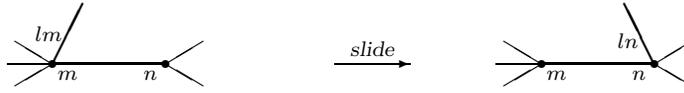
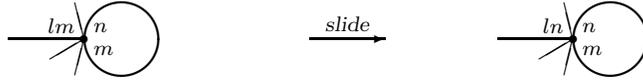
\begin{pict}{100}{11}
\thicklines
\put(75,3){\circle*{1}}
\put(90,3){\circle*{1}}
\put(75,3){\line(1,0){15}}
\put(90,3){\line(-1,2){4}}

\thinlines
\put(47.5,3){\vector(1,0){10}}
\zindex{52.5}{5}{\mbox{slide}}

\put(90,3){\line(5,3){5}}
\put(90,3){\line(5,-3){5}}
\put(75,3){\line(-5,3){5}}
\put(69,3){\line(1,0){6}}
\put(75,3){\line(-5,-3){5}}

\scriptsize
\zindex{77}{1.5}{m}
\zindex{88}{1.5}{n}
\zindex{86.5}{6}{ln}

\thicklines
\put(10,3){\circle*{1}}
\put(25,3){\circle*{1}}
\put(10,3){\line(1,0){15}}
\put(10,3){\line(1,2){4}}

\thinlines
\put(25,3){\line(5,3){5}}
\put(25,3){\line(5,-3){5}}
\put(10,3){\line(-5,3){5}}
\put(4,3){\line(1,0){6}}
\put(10,3){\line(-5,-3){5}}

\scriptsize
\zindex{12}{1.5}{m}
\zindex{23}{1.5}{n}
\zindex{9.5}{7}{lm}
 
\end{pict} 
or

\begin{pict}{100}{10}
\thicklines

\put(82.5,5){\circle*{1}}
\put(87.5,5){\circle{10}}
\put(72.5,5){\line(1,0){10}}

\thinlines
\put(47.5,5){\vector(1,0){10}}
\zindex{52.5}{7}{\mbox{slide}}

\put(82.5,5){\line(-5,-3){4.5}}
\put(82.5,5){\line(-1,-4){1.2}}
\put(82.5,5){\line(-1,4){1.2}}

\scriptsize
\zindex{85.1}{3.5}{m}
\zindex{84.7}{6.5}{n}
\zindex{79.9}{6.8}{ln}

\thicklines

\put(17.5,5){\circle*{1}}
\put(22.5,5){\circle{10}}
\put(7.5,5){\line(1,0){10}}

\thinlines
\put(17.5,5){\line(-5,-3){4.5}}
\put(17.5,5){\line(-1,-4){1.2}}
\put(17.5,5){\line(-1,4){1.2}}

\scriptsize
\zindex{20.1}{3.5}{m}
\zindex{19.7}{6.5}{n}
\zindex{14.5}{6.8}{lm}

\end{pict} 
\end{proposition}

The proof is straightforward and is left to the reader. 

\begin{example} \label{eg} 
Choose non-zero integers $m$, $n$, $r$, and $s$. 
The following diagrams depict a sequence of elementary deformations
between generalized Baumslag--Solitar trees. The initial and final trees
are the examples that interest us. Call them $X$ and $Y$
respectively, and let $G$ be the group. 

\medskip

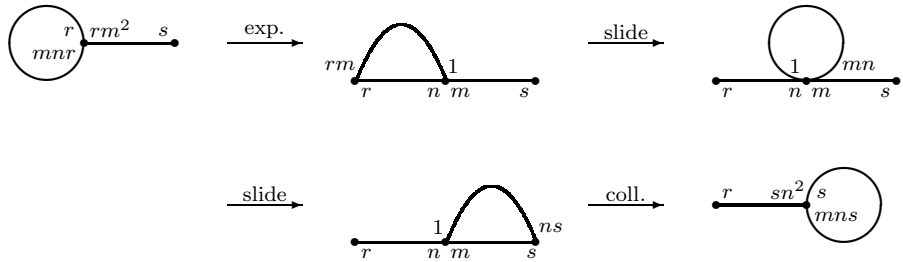
\begin{pict}{120}{13}
\thicklines
\put(10,8){\circle*{1}}
\put(22,8){\circle*{1}}
\put(5,8){\circle{10}}

\put(10,8){\line(1,0){12}}
\scriptsize
\zindex{6}{6.5}{mnr}
\zindex{8}{9.5}{r}
\zindex{13.5}{10}{rm^2}
\zindex{20.5}{9.5}{s}

\put(46,3){\circle*{1}}
\put(58,3){\circle*{1}}
\put(70,3){\circle*{1}}

\put(46,3){\line(1,0){24}}
\qbezier(46,3)(52,18)(58,3)
\zindex{47.5}{1.5}{r}
\zindex{56.5}{1.5}{n}
\zindex{60}{1.5}{m}
\zindex{68.5}{1.5}{s}
\zindex{44}{5}{rm}
\zindex{59}{5}{1}
 
\put(94,3){\circle*{1}}
\put(106,3){\circle*{1}}
\put(118,3){\circle*{1}}

\put(94,3){\line(1,0){24}}
\put(106,8){\circle{10}}
\zindex{95.5}{1.5}{r}
\zindex{104.5}{1.5}{n}
\zindex{108}{1.5}{m}
\zindex{116.5}{1.5}{s}
\zindex{104.5}{5}{1}
\zindex{113}{5}{mn}

\thinlines
\put(29,8){\vector(1,0){10}}
\zindex{34}{9.5}{\mbox{exp.}}

\put(77,8){\vector(1,0){10}}
\zindex{82}{9.5}{\mbox{slide}}

\end{pict}%

\smallskip

\begin{pict}{120}{13}
\thicklines
\put(94,8){\circle*{1}}
\put(106,8){\circle*{1}}
\put(111,8){\circle{10}}

\put(94,8){\line(1,0){12}}
\scriptsize
\zindex{110}{6.5}{mns}
\zindex{108}{9.5}{s}
\zindex{103.5}{10}{sn^2}
\zindex{95.5}{9.5}{r}

\put(46,3){\circle*{1}}
\put(58,3){\circle*{1}}
\put(70,3){\circle*{1}}

\put(46,3){\line(1,0){24}}
\qbezier(58,3)(64,18)(70,3)
\zindex{47.5}{1.5}{r}
\zindex{56.5}{1.5}{n}
\zindex{60}{1.5}{m}
\zindex{69.5}{1.5}{s}
\zindex{57}{5}{1}
\zindex{72}{5}{ns}

\thinlines
\put(29,8){\vector(1,0){10}}
\zindex{34}{9.5}{\mbox{slide}}

\put(77,8){\vector(1,0){10}}
\zindex{82}{9.5}{\mbox{coll.}}

\end{pict}%
\end{example}

\begin{proposition} \label{slide} 
If $m \nmid n$ and $n \nmid m$ then $X$ and $Y$ are not related
by slide moves. 
\end{proposition}

\begin{proof}
Note that $m \not= \pm 1$, and so $rm^2 \nmid mnr$, $mnr \nmid rm^2$, and
$rm^2 \nmid r$.  Thus $X$ admits only one slide move, in which the free
edge travels around the loop. This move changes the index $rm^2$ to $rm^3
n$. Repeating this move, the index becomes $rm^{k+2}n^k$ with $k \geq
0$. Call this $G$--tree $X_k$. Then since $rm^{k+2}n^k \nmid mnr$ and
$rm^{k+2}n^k \nmid r$, the only slide moves available from $X_k$ are those
resulting in $X_{k+1}$ and $X_{k-1}$ (when $k \geq 1$). Since $Y \not\cong
X_k$ for any $k$, the result follows. 
\end{proof}

\begin{remark}
If one considers the set of reduced $G$--trees for a given group $G$,
then the example shows that the relation of slide--equivalence may be
strictly finer than the relation of elementary deformation. 
\end{remark}

Throughout the rest of this section we verify that $X$ and $Y$ are in
fact JSJ decompositions of $G$. This will be the case as long as $r,s
\not= \pm 1$. First we review some basic properties of generalized
Baumslag--Solitar groups. The key property, from the point of view of JSJ
decompositions, is given in Lemma \ref{gbssplittings}. 

\begin{definition}
Let $X$ be a $G$--tree. An element $\gamma \in G$ is \emph{elliptic} if
it fixes a vertex of $X$ and \emph{hyperbolic} otherwise. If $\gamma$ is
hyperbolic then there is a unique $\gamma$--invariant line in $X$, called
the \emph{axis} of $\gamma$, on which $\gamma$ acts as a translation
\cite[Chapter I, Proposition 24]{serre:trees}. From this description it
is clear that for any $n\not= 0$, the element $\gamma^n$ is hyperbolic if
and only if $\gamma$ is, and when this occurs they have the same axis. 

Two elements $\gamma,\delta \in G$ are \emph{commensurable} if there
exist non-zero integers $m,n$ such that $\gamma^m = \delta^n$. Note that
commensurable hyperbolic elements have the same axis. The
\emph{commensurator} of $\gamma$ is the set of all elements $\delta \in
G$ such that $\delta \gamma \delta^{-1}$ and $\gamma$ are commensurable. 
\end{definition}

\begin{lemma} \label{commensurator} 
Let $X$ be a $G$--tree. If $\gamma \in G$ is hyperbolic then its
commensurator stabilizes its axis. If $\gamma \in G$ is elliptic and $X$
is a generalized Baumslag--Solitar tree, then the commensurator of
$\gamma$ is all of $\, G$. 
\end{lemma}

\begin{proof}
Suppose that $\delta \gamma \delta^{-1}$ and $\gamma$ are commensurable,
where $\gamma$ is hyperbolic
with axis $L$. Then $\delta \gamma \delta^{-1}$ also has axis
$L$. However, the axis of $\delta \gamma \delta^{-1}$ is $\delta L$, and
hence $\delta$ stabilizes $L$.  Next suppose that $X$ is a generalized
Baumslag--Solitar tree. Then all non-trivial elliptic elements are
commensurable, and hence every elliptic element is commensurable with all
of its conjugates. 
\end{proof}

\begin{lemma} \label{gbsproperties} 
Let $X$ be a generalized Baumslag--Solitar tree with group $G \not\cong
\Z$. 
Then: 
\begin{enumerate}
\item \label{p1.5} $G$ is not free; 
\item \label{p2} $G$ is torsion--free and has cohomological dimension $2$; 
\item \label{p1} $G$ has one end, if it is finitely generated;
\item \label{p3} $X$ contains a $G$--invariant line if and only if $G$ is
isomorphic to $\Z \times \Z$ or the Klein bottle group. 
\end{enumerate}
\end{lemma}

\begin{proof}
For (\ref{p1.5}), first suppose that any two vertex stabilizers are
contained in an infinite cyclic subgroup.  If there is a maximal
stabilizer $C$, then $G \cong C \times \pi_1(G \backslash X)$, which is
either $\Z$ or non-free. If there is no maximal stabilizer, then the
set of elliptic elements is an abelian, non-finitely generated subgroup
of $G$. Free groups contain no such subgroups. Finally, if there are two
stabilizers not contained in an infinite cyclic subgroup, let $\gamma$
and $\delta$ be their generators, respectively. Then $\gamma$ and
$\delta$ cannot generate a free subgroup of $G$ because they are
commensurable. 

Claim (\ref{p2}) follows from the fact that $G$ is the fundamental group
of a space $Z$ with universal cover homeomorphic to $X \times \R$. The
space $Z$ is the total space of a graph of spaces in which every 
vertex and edge space is a circle (see \cite{scottwall} for these
notions). Note that $Z$ is aspherical and $2$--dimensional. 

For (\ref{p1}), if $G$ is finitely generated then it acts cocompactly on
the minimal subtree of $X$. Replace $X$ by this subtree. Then the space
$Z$ is compact, and so $G$ and $X \times \R$ have the same number of
ends. The space $X \times \R$ has one end unless $X$ is compact, but this
does not occur because $G \not\cong \Z$. 

For (\ref{p3}), if $X$ has an invariant line then $G$ acts on this line
with infinite cyclic stabilizers. Since $G \not\cong \Z$, this
action is cocompact. Hence $G = \Z \ast_{\Z} \Z$ or $G =
\Z \, \ast_{\Z}$, with inclusion maps multiplication by $\pm 2$
in the first case and $\pm 1$ in the second. For the converse, if $G
\cong \Z \times \Z$ or the Klein bottle group, then no minimal
subtree can have more than two ends, for this would imply that $G$ has
exponential growth.
\end{proof}

\begin{lemma} \label{gbssubgroups} 
Let $X$ be a generalized Baumslag--Solitar tree with group $G$. Every
subgroup $H$ of $\, G$ is either a generalized Baumslag--Solitar group or a
free group (and not both except for $\Z$). If $H$ is free and
non-abelian then every non-trivial element of $H$ is hyperbolic. 
\end{lemma}

\begin{proof}
Every edge stabilizer has finite index in its neighboring vertex
stabilizers, and hence $X$ is locally finite. This implies that for every
group action on $X$, all vertex and edge stabilizers are
commensurable. Thus the stabilizers of $H$ acting on $X$ are either all
infinite cyclic or all trivial. If the former occurs then $H$ is not free
unless it is $\Z$, by Lemma \ref{gbsproperties}(\ref{p1.5}). Hence
if $H$ is free and non-abelian then its stabilizers in $X$ are all
trivial. 
\end{proof}

\begin{definition}
A \emph{splitting of $\, G$ over $\, C$} is a graph of groups
decomposition of the form $G = A \, \ast_C B$ with $A \not= C \not= B$,
or $G = A \, \ast_C$ (with no restriction on $C$). We also say that $G$
\emph{splits} over $C$. If $\As_1$ is a splitting of $G$ over $C$ and
$\As_2$ is another graph of groups decomposition of $G$, we say that
$\As_1$ is \emph{hyperbolic in $\As_2$} if some element of $C$ is
hyperbolic in the action of $G$ on the Bass--Serre tree of $\As_2$. 
\end{definition}

\begin{lemma} \label{gbssplittings} 
Let $X$ be a generalized Baumslag--Solitar tree with group $G$. Assume
that $X$ contains no $G$--invariant point or line. If $G$ splits over a
$2$--ended group $C$, then $C$ is contained in a vertex stabilizer of $\,
X$. 
\end{lemma}

\begin{proof}
Let $T$ be the Bass--Serre tree for the splitting over $C$, and let $e$
be an edge of $T$ with stabilizer $C$. Note that $T$ contains no
invariant point or line, for this would imply that $G \cong \Z$, or
that $G$ splits as $A_1 \, \ast_C B_1$ with $[A_1:C] = [B_1:C] = 2$, or
$A_1 \, \ast_C$ with both inclusions surjective. In these cases $G = \Z
\times \Z$ or the Klein bottle group, contradicting the 
hypotheses (by Lemma \ref{gbsproperties}(\ref{p3})). 

Let $\gamma \in G$ be a generator of a vertex stabilizer in $X$. By Lemma
\ref{commensurator} the commensurator of $\gamma$ is all of $G$. Now
consider the action of $G$ on $T$. If $\gamma$ is hyperbolic then its
commensurator stabilizes its axis, but we have observed that $G$ does not
stabilize any line in $T$. Hence $\gamma$ fixes a vertex $v$ of $T$.  Now
let $\delta \in G$ be chosen so that $e$ separates $v$ from $\delta
v$. Since $\gamma$ and $\delta \gamma \delta^{-1}$ are commensurable,
some power of $\gamma$ fixes both $v$ and $\delta v$, and also $e$. Thus
$\gamma^n \in C$ for some $n$. Note that $C$ is infinite cyclic because
it is torsion free and $2$--ended. Letting $c$ be a generator, we have
that $c^m = \gamma^n$ for some $m$. Since this element is elliptic in
$X$, $c$ is also elliptic in $X$. 
\end{proof}

\begin{definition}
In order to discuss JSJ decompositions we need a further definition, from
\cite{ripssela}. A subgroup $H$ of $\, G$ is \emph{quadratically hanging}
if $H$ is isomorphic to the fundamental group of a compact $2$--orbifold
with boundary, and there exists a minimal $G$--tree with the following
properties: all edge stabilizers are infinite cyclic, there is a vertex
$v$ with stabilizer $H$, and the stabilizers of the edges incident to $v$
are precisely the conjugates (in $H$) of the peripheral subgroups of
$H$. It is also required that the $2$--orbifold have negative Euler
characteristic and contain a pair of intersecting ``weakly essential''
simple closed curves. 
\end{definition}

\begin{lemma} \label{nocmq} 
Let $G$ be a generalized Baumslag--Solitar group. Then $G$ contains no
quadratically hanging subgroups. 
\end{lemma}

\begin{proof}
Let $X$ be a generalized Baumslag--Solitar tree with group $G$.  Suppose
that $H \subseteq G$ is a quadratically hanging subgroup with
corresponding $G$--tree $T$.  Since $G$ is torsion--free, $H$ is the
fundamental group of a compact surface with boundary and negative Euler
characteristic. 
It is conceivable that the surface is closed, in which case $T$ is a
point and $H = G$. However, no generalized Baumslag--Solitar group is
isomorphic to a closed surface group of this type. 
For example, no generalized Baumslag--Solitar group is hyperbolic,
except for $\Z$. 

Thus, we can assume that the surface has non-empty boundary, so that $H$
is a non-abelian free group. In particular $G$ is not isomorphic to $\Z$,
$\Z \times \Z$, or the Klein bottle group, and so $X$ 
contains no $G$--invariant point or line. Now let $C\subseteq H$ be a
peripheral subgroup, and consider the splitting of $G$ over $C$ arising
from the $G$--tree $T$. By Lemma \ref{gbssplittings}, every element of
$C$ is elliptic relative to $X$. This contradicts Lemma
\ref{gbssubgroups}. 
\end{proof}

\begin{definition}
The next theorem \cite[Theorem 7.1]{ripssela} defines the notion of a JSJ
decomposition. To be specific, a \emph{JSJ decomposition} of a group is a
graph of groups decomposition satisfying the properties of $\, \Gs$ in
Theorem \ref{rsthm}.  The precise definition of \emph{CMQ subgroup}
is not important here, except that it is a quadratically hanging subgroup
with additional properties. Similarly, we omit the definition of
\emph{weakly essential} simple closed curve. 

A splitting $\As_1$ of $\, G$ is \emph{unfolded} if there do not exist a
splitting $\As_2$ and a non-trivial fold $\As_2 \to \As_1$. Note that
such a fold involves only $G$--trees having a single edge orbit, and
therefore induces an isomorphism of quotient graphs. 
A graph of groups $\Gs_1$ is \emph{unfolded} if every
splitting arising from its edges is unfolded in the above sense. Note
that this is very different from requiring no non-trivial fold $\Gs_2 \to
\Gs_1$ to exist. 
\end{definition}

\begin{theorem}[Rips--Sela] \label{rsthm} 
Let $G$ be a finitely presented group with one end. There exists a
reduced, unfolded graph of groups decomposition $\, \Gs$ of $\, G$ with
infinite cyclic edge groups, such that the following conditions hold. 
\begin{enumerate}

\item \label{rs1} Every canonical maximal quadratically hanging (CMQ)
subgroup of $\, G$ is conjugate to a vertex group of $\, \Gs$. Every
quadratically hanging subgroup of $\, G$ can be 
conjugated into one of the CMQ subgroups of $\, G$. Every non-CMQ vertex
group of $\, \Gs$ is elliptic in every graph of groups decomposition of
$\, G$ having infinite cyclic edge groups. 

\item \label{rs2} A splitting of $\, G$ over $\Z$ which is hyperbolic
in another splitting over $\Z$ is obtained from $\, \Gs$ by cutting a
$2$--orbifold corresponding to a CMQ subgroup of $\, G$ along a weakly
essential simple closed curve. 

\item \label{rs3} Let $\As$ be a splitting of $\, G$ over $\Z$ which
is elliptic with respect to every other splitting over $\Z$. Then
there exists a $G$--equivariant simplicial map from a subdivision of
$T_{\Gs}$, the Bass--Serre tree of $\, \Gs$, to $T_{\As}$ (the Bass--Serre
tree of $\As$). 

\item \label{rs4} Let $\As$ be a graph of groups decomposition of $\, G$
with infinite cyclic edge groups. Then there exist a decomposition
$\Gs'$ obtained from $\Gs$ by splitting the CMQ subgroups along weakly
essential simple closed curves on their corresponding $2$--orbifolds, and
a $G$--equivariant simplicial map from a subdivision of $\, T_{\Gs'}$ to
$T_{\As}$. 

\end{enumerate}
\end{theorem}
Rips and Sela also include a uniqueness statement in their theorem, which
we discuss in Section \ref{uniquesec}. 

\begin{remark}
Our notion of ``reduced'' $G$--trees is stronger than the definition used
by Rips and Sela, and stated in \cite{bestvina:accessibility}. However it
is true that if one performs collapse moves on a JSJ decomposition then
the result is again a JSJ decomposition. Thus Theorem \ref{rsthm} is
valid with our definition. This issue will be relevant when we discuss
uniqueness in the next section. 
\end{remark}

\begin{theorem} \label{gbsjsj} 
Let $X$ be a generalized Baumslag--Solitar tree with group $G$. If $X$ is
reduced, unfolded, and is not a point or line, then $X$ is a JSJ
decomposition of $\, G$. 
\end{theorem}

\begin{proof}
Reduced trees are minimal, so $X$ contains no invariant point or
line. The first two statements of (\ref{rs1}) are vacuously true, by
Lemma \ref{nocmq}. The third statement is a consequence of Lemma
\ref{gbssplittings}, as follows. Let $\As$ be a graph of groups
decomposition of $\, G$ with infinite cyclic edge groups. Choose an edge
group $C \subseteq G$ of $\As$. Then $C$ is elliptic in $X$ by Lemma
\ref{gbssplittings}. As $C$ is infinite cyclic, each vertex stabilizer of
$X$ is commensurable with $C$. Therefore each vertex stabilizer of $X$ is
elliptic in $\As$.  Claim (\ref{rs2}) is also vacuously true by Lemma
\ref{gbssplittings}; if $G$ splits over infinite cyclic subgroups $C_1$
and $C_2$, then the lemma implies that $C_1$ and $C_2$ are commensurable,
and hence neither splitting is hyperbolic in the other.  Claim
(\ref{rs3}) is a special case of (\ref{rs4}), since there are no CMQ
subgroups. To prove (\ref{rs4}) it suffices to verify that every vertex
stabilizer of $X$ is elliptic in $\As$. This was just shown in the proof
of (\ref{rs1}). 
\end{proof}
	
\begin{remark}
A similar argument shows that any generalized Baumslag--Solitar tree as
above is also a JSJ decomposition over $2$--ended groups in the sense of
Dunwoody and Sageev \cite{dunwoodysageev}, and a JSJ decomposition over
slender groups in the sense of Fujiwara and Papasoglu
\cite{fujiwara:jsj}. 
\end{remark}

\begin{proposition} \label{gbsunfolded} 
Let $X$ be a cocompact generalized Baumslag--Solitar tree with group
$G$. If every edge stabilizer is a proper subgroup of its neighboring
vertex stabilizers then $X$ is unfolded. 
\end{proposition}

The assumption of cocompactness is not necessary, but it allows for a
simpler proof. 

\begin{proof}
Let $e$ be an edge of $X$ and let $\widehat{X}$ be the tree obtained from
$X$ by collapsing each connected component of $X - G e$ to a vertex. Call
this quotient map $q \co X \to \widehat{X}$. Then $\widehat{X}$ is the
Bass--Serre tree corresponding to the splitting of $\, G$ associated to
$e$. Now suppose that there is a non-trivial fold $f \co Y \to
\widehat{X}$. Let $e' \in E(Y)$ be an edge with $f(e') = q(e)$.  Without
loss of generality, assume that the fold occurs at $\partial_0 e'$. 

Consider the stabilizer of the vertex $\partial_1 e$ in $X$, and let
$\gamma$ be a generator. Then $G_e$ is generated by $\gamma^n$ for some
$n$. Also $G_{e'} \subsetneq G_{q(e)} = G_e$ and $G_{e'} \not= \{1\}$, as
$G$ has one end (by cocompactness and Lemma
\ref{gbsproperties}(\ref{p1})). Thus $G_{e'}$ is 
generated by $\gamma^{mn}$ for some 
$m>1$. This implies that $\gamma^n e' \not= e'$, and $\partial_0 \gamma^n
e' = \partial_0 e'$ because the fold occurs at $\partial_0 e'$. The edges
$\gamma^n e'$ and $e'$ separate $Y$ into three connected components, 
adjacent to the vertices $\partial_0 e'$, $\partial_1 e'$, and
$\partial_1 \gamma^n e'$. Call these subtrees $Y_0$, $Y_1$ and
$\gamma^n(Y_1)$ accordingly. It is clear that $\gamma$ fixes no point of
$Y_1$ or $\gamma^n(Y_1)$. Also define the subtrees $\widehat{X}_0$,
$\widehat{X}_1$, $X_0$, and $X_1$ similarly, as the connected components of
$\widehat{X} - q(e)$ and $X - e$. 

We claim that $\gamma e = e$, and hence $G_e = G_{\partial_1 e}$. Note
that since $\gamma^{mn}$ is elliptic in $Y$, so is $\gamma$. Therefore
$\gamma$ fixes a vertex in $Y_0$. By equivariance of $f$ it also fixes a
vertex $v$ of $\widehat{X}_0$. Then $\gamma$ stabilizes the subtree
$q^{-1}(v) \subseteq X_0$. Since $e$ separates $q^{-1}(v)$ from
$\partial_1 e$ and $\gamma$ stabilizes both, $\gamma$ also fixes $e$. 
\end{proof}

Finally we return to the main example. If $r,s \not= \pm 1$ then $X$ and
$Y$ are reduced and unfolded. Since neither tree is a point or
line, both are JSJ decompositions of $\, G$.

\section{Uniqueness of decompositions} \label{uniquesec} 

We have seen that JSJ decompositions of a given group need not be related
by slide moves. However, the results of \cite{forester:trees} imply
that they are unique up to elementary deformation, and also that
many JSJ decompositions are genuinely unique. The reasoning applies to
other decompositions as well: the JSJ decompositions of Dunwoody and
Sageev \cite{dunwoodysageev} and Fujiwara and Papasoglu
\cite{fujiwara:jsj}, and also accessible (or one-ended) decompositions. 

The other JSJ decompositions just mentioned are defined in a similar way
to the Rips--Sela version, by means of a kind of universal property. In
Theorem \ref{rsthm} this property refers to splittings over infinite cyclic
groups. The JSJ decomposition has edge groups of this type, and it
provides a simultaneous description of every such splitting. In a similar
fashion, the JSJ decomposition of Fujiwara and Papasoglu deals with
splittings over slender groups. (A group is \emph{slender} if every
subgroup is finitely generated.) For the Dunwoody--Sageev decomposition
one must first choose a ``closed class'' $\Cs$ of slender groups, and
then the JSJ decomposition refers to splittings over elements of
$\Cs$. For example, $\Cs$ could be the class of $2$--ended groups, or the
class of finite extensions of $\Z \times \Z$. For further details on
these JSJ decompositions see \cite{dunwoodysageev} and
\cite{fujiwara:jsj}. 

As mentioned earlier we require JSJ decompositions to be reduced, in
order to apply Theorem \ref{rigidity} below. In all
three versions this is easily arranged by performing collapse moves,
though in the case of the Dunwoody--Sageev decomposition a small
modification is required. Namely, we must use their decomposition
$\Gs_{\textup{red}}$ rather than $\Gs$, which is obtained from
$\Gs_{\textup{red}}$ by subdivision.  

An \emph{accessible decomposition} is a reduced graph of groups whose edge
groups are finite and whose vertex groups each have at most one
end. Dunwoody showed in \cite{dunwoody:accessibility} that every finitely
presented group has an accessible decomposition. 

The property shared by all of these decompositions is that any two
particular decompositions (of the same kind) have the same elliptic
subgroups. Here, an \emph{elliptic subgroup} is any subgroup that fixes a
vertex of the given $G$--tree. To state our uniqueness theorem we need
one additional notion. A $G$--tree is \emph{strongly slide--free} if it
is minimal and, for all edges $e$ and $f$ with $\partial_0 e = \partial_0
f$, $G_e \subseteq G_f$ implies $f \in Ge$. In terms of graphs of groups
this means that for every vertex group $A$, if $C$ and $C'$ are
neighboring edge groups then no conjugate (in $A$) of $C$ is contained in
$C'$. Many graphs of groups arising in nature have this property, though
there are obvious exceptions such as ascending HNN extensions. 

\begin{theorem} \label{uniqueness} 
Let $G$ be a finitely generated group. Suppose that $X$ and $Y$ are
$G$--trees representing decompositions of one of the following types: JSJ
decompositions in the sense of Rips and Sela, Dunwoody and Sageev, or
Fujiwara and Papasoglu, or accessible decompositions. Then $X$ and $Y$
are related by an elementary deformation. If $X$ is strongly slide--free
then there is a unique $G$--isomorphism $X \to Y$. 
\end{theorem}

Thus, in each case, strongly slide--free decompositions are genuinely
unique. The proof is a direct application of the following two results
\cite[Theorems 1.1 and 1.2]{forester:trees}. 

\begin{theorem}
Let $G$ be a group and let $X$ and $Y$ be cocompact $G$--trees. Then $X$
and $Y$ are related by an elementary deformation if and only if they have
the same elliptic subgroups. \endproof 
\end{theorem}

\begin{theorem} \label{rigidity} 
Let $X$ and $Y$ be cocompact $G$--trees that are related by an elementary
deformation. If $X$ is strongly slide--free and $Y$ is reduced then
there is a unique $G$--isomorphism $X \to Y$. \endproof 
\end{theorem}

\begin{proof}[Proof of Theorem \ref{uniqueness}] 
Both trees are cocompact because $G$ is finitely generated. Now it
suffices to verify that the elliptic subgroups for $X$ and $Y$ agree. An
equivalent property is that there exist equivariant maps between $X$ and
$Y$ in each direction. This is proved in \cite[Theorem 7.1(v)]{ripssela}
and \cite[p. 43]{dunwoodysageev} respectively for the first two types of
decomposition. The third case is similar to these, and it can be derived
formally from the main theorem of \cite{fujiwara:jsj}. 

For accessible decompositions we argue as follows. Let $G_x$ be a vertex
stabilizer in $X$ and consider its action on $Y$. This action has finite
edge stabilizers, and since $G_x$ has at most one end, the action must be
trivial. Thus, every elliptic subgroup in $X$ fixes a vertex in $Y$, and
conversely by symmetry. 
\end{proof}

\providecommand{\bysame}{\leavevmode\hbox to3em{\hrulefill}\thinspace}

\end{document}